\documentclass[11pt]{article}
\title{On invariant subspaces of a linear operator}
\author{M.\,I. Belishev and S.\,A. Simonov}
\date{}

\usepackage[a5paper, margin={1cm, 2cm}]{geometry}

\usepackage[cp1251]{inputenc}
\usepackage{amssymb,amsmath,amsthm}
\usepackage{ulem}
\usepackage{graphicx,color}
\usepackage{mathrsfs}
\usepackage{array, amsfonts, mathrsfs}
\usepackage{epstopdf}
\usepackage[colorlinks=true]{hyperref}
\usepackage{float}
\usepackage{enumitem}
\usepackage{amsmath}
\usepackage[dvipsnames]{xcolor}

\newtheorem{Theorem}{Theorem}

\newtheorem{Definition}{Definition}

\def\mG{\mathscr G}
\def\mH{\mathscr H}

\def\mL{\mathscr L}

\def\Dom{{\rm Dom\,}}

\def\harp{\upharpoonright}

\def\bul{\noindent$\bullet$\,\,\,}

\begin{document}
\maketitle

\begin{abstract}
	We discuss the concept of invariant subspaces for unbounded linear operators, point out some shortcomings of known definitions, and propose our own.
\end{abstract}

\subsubsection*{About the work}

This note discusses the concept of an invariant subspace for an {\it unbounded} linear operator. We highlight some drawbacks of known definitions \cite{AG_engl,BirSol,DSch,Ed,Kato} and propose our own, which, in our opinion, avoids these shortcomings.

Citations are given in notation that does not always match the original one: this is done for readability. 

\subsubsection*{Bounded operators}

Let $\mH$ be a (complex and separable) Hilbert space, $\mG\subset\mH$ a (closed) subspace, $P_\mG$ the orthogonal projector in $\mH$ onto $\mG$, and $\mG^\bot:=\mH\ominus\mG$. The following definitions and facts are standard: see the fundamental monographs \cite{DSch,KAk,BirSol}.

\begin{Definition}\label{Def 1}
	A subspace $\mG\not=\{0\}$ is invariant for a bounded operator $A$ if the relation $A\mG\subset\mG$ holds, or equivalently, the equality $AP_\mG=P_\mG AP_\mG$. In this case, the operator $A_\mG:\mG\to\mG$, $A_\mG x:=Ax$, is called the part of $A$ in the subspace $\mG$. If both $\mG$ and $\mG^\bot$ are invariant, we say that $\mG$ (and $\mG^\bot$) reduces $A$, which is equivalent to the commutation $AP_\mG=P_\mG A$.
\end{Definition}

If ${\rm dim\,}\mG<\infty$, then the part $A_\mG$ has eigenvectors in $\mG$, satisfying $Ax=\lambda x$, $x\not=0$, $\lambda\in\Bbb C$.

\subsubsection*{Unbounded operators}

\bul The most general definition is given in \cite{AG}, Chapter IV, §45:
\begin{Definition}\label{Def AG}
	A subspace $\mG$ is invariant for $A$ if the following relation holds:
	\begin{equation*}
		A \ [\mG\cap\Dom A]\,\subset\mG.
	\end{equation*}
	In this case, the operator $A_\mG:\mG\to\mG$, $A_\mG x:=Ax$, $x\in \mG\cap\Dom A$, is called the part of $A$ in $\mG$.
\end{Definition}

Note that in the translation \cite{AG_engl}, the term ``part'' is replaced with ``restriction'', which is hardly fortunate: ``restriction'' usually refers to the operator $A\harp\mG:\mG \to \mH$.

The above definition is not without drawbacks. Indeed, let an element $x\in\Dom A$ satisfy $Ax\not\in\Dom A$ and $x_1, x_2,\dots, x_n\not\in\Dom A$. Then the subspace ${\rm span\,}\{x,Ax, x_1,\dots,x_n\}$ turns out to be invariant for $A$, and thus the operator acquires many ``unnecessary'' invariant subspaces. Another example is ${\rm span\,}\{x, x_1,\dots,x_n\}$ with $x\not= 0$ under the condition $Ax=\lambda x$, $\lambda\in\Bbb C$.

A similar argument is given in \cite{Kato}, Chapter III, §5: ``It is rather difficult to extend the definition of an invariant subspace to unbounded operators in $\mH$, since the inclusion $A\mG\subset\mG$ will hold whenever $\mG$ contains only the zero vector from $\Dom A$''.

We add that in the Russian original \cite{AG}, below the definition of an invariant subspace (p. 131), it is stated: ``Note that every invariant subspace of finite dimension contains at least one eigenvector, as we know from linear algebra''. The example ${\rm span\,}\{x,Ax, x_1,\dots,x_n\}$ shows that this is obviously false. This statement is absent in the English translation.
\smallskip

\bul In \cite{Ples}, Chapter VI, §7, invariance is defined as follows.
\begin{Definition}\label{Def Ples}
	A subspace $\mG$ is invariant for $A$ if the following relations hold:
	\begin{equation}\label{Eq def Ples}
		1)\,\,\,A \ [\mG\cap\Dom A]\,\subset\mG, \qquad 2)\,\,\,P_\mG\Dom A\subset\Dom A.
	\end{equation}
	If only 1) holds, we say that $\mG$ is weakly invariant. A subspace $\mG$ (orthogonally) reduces $A$ if both $\mG$ and $\mG^\bot$ are invariant for $A$. In the latter case, the operator $A_\mG$ induced by $A$ in $\mG$ is called the part of $A$.
\end{Definition}

Thus, {\it weak invariance} is invariance according to Definition \ref{Def AG}. Reducibility means that when projecting $\Dom A$ onto $\mG$, the elements of $\Dom A$ remain in $\Dom A$, with $A_{\mG^\bot}$ also being a part of $A$ (in $\mG^\bot$), and we write $A=A_{\mG}\oplus A_{\mG^\bot}$.

In \cite{BirSol}, invariance itself is not defined, but reducibility is understood as in Definition \ref{Def Ples}. However, in our opinion, the condition $P_\mG\Dom A\subset\Dom A$ is too restrictive for applications. Let us illustrate this with an example.
\smallskip

Let $\mH=L_2(\Bbb R)$, $\tilde L=\sum_{k=0}^n p_k(\cdot)D^k$ be a differential expression with locally bounded coefficients, defining the differential operator $L:\mH\to\mH$, $\Dom L= W_2^n(\Bbb R)$, $Ly:=\tilde L y$. The subspaces $\mH^{ab}:=\{y\in\mH\,|\,\, {\rm supp\,}y\subset[a,b]\},\,\, -\infty\leqslant a<b\leqslant\infty$, are invariant according to Definition \ref{Def AG} but not according to Definition \ref{Def Ples}, since the second condition in (\ref{Eq def Ples}) fails if $a$ or $b$ are finite.

\subsubsection*{Our proposals}

\bul The following formulations seem a reasonable compromise, motivated by the content and terminology of the monograph \cite{Naimark}.
\begin{Definition}\label{Def Our}
	A subspace $\mG$ is invariant for $A$ if the following relations hold:
	\begin{equation*}
		1)\,\,\,\overline{\mG\cap\Dom A}=\mG, \qquad 2)\,\,\,A\ [\mG\cap\Dom A]\subset\mG.
	\end{equation*}
	In this case, the operator $A_\mG:\mG\to\mG$, $\Dom A_\mG={\mG\cap\Dom A}$, $A_\mG y:=Ay$, is called the part of $A$ in $\mG$. If both $\mG$ and $\mG^\bot$ are invariant, we say that $\mG$ (and $\mG^\bot$) splits $A$ into parts $A_\mG$ and $A_{\mG^\bot}$ and write $A\supset A_\mG\oplus A_{\mG^\bot}$. If, in addition, $P_\mG\Dom A\subset \Dom A$ holds, we say that $\mG$ (and $\mG^\bot$) reduces $A$ and write $A=A_\mG\oplus A_{\mG^\bot}$.
\end{Definition}
\noindent The following facts and results confirm the adequacy of the proposed definitions.
\smallskip

$\star$\,\,\, If $A$ is bounded, then Definitions \ref{Def 1} and \ref{Def Our} give the same notion of invariance. In this case, the relations $A\supset A_\mG\oplus A_{\mG^\bot}$ and $A=A_\mG\oplus A_{\mG^\bot}$ are equivalent. If $A$ is a closed operator, then $A_\mG$ is also closed.
\smallskip

$\star$\,\,\, The subspaces $\mH^{ab}$ are invariant and split the differential operator $L$ but do not reduce it.
\smallskip

$\star$\,\,\, The following result appears quite natural. All terms are understood in the sense of Definition \ref{Def Our}.
\begin{Theorem}
	If the operator $A$ is symmetric (i.e., $A\subset A^*$), $\mG$ is its invariant subspace, and its part in $\mG$ is self-adjoint (i.e., $A_\mG=A_\mG^*$ in $\mG$), then $\mG$ reduces $A$.
\end{Theorem}
This statement can be derived directly from our definitions or extracted from more general considerations in \cite{Behrndt}.
\smallskip

\bul There is a specific situation where another operator associated with $A$ can be considered as its part. Let a linear manifold $\mL\subset\mG$ be such that $\overline{\mL}=\mG$, $\mL\subset\Dom A$, and $A\mL\subset\mL$. It seems natural to regard the operator $A\harp\mL$ as a part of $A$ on the linear manifold $\mL$. At the same time, $A$ may or may not have a part $A_\mG$ (in the sense of Definition \ref{Def Our}), and in the first case, we obviously have $A_\mL\subset A_\mG$.

An example is the so-called {\it wave part} of a symmetric operator, where a class of solutions to a dynamical system ({\it smooth waves}) plays the role of $\mL$; see \cite{B_DSBC_3,BSim_charact_SL}.
\smallskip

\bul The following concluding proposal summarizes our considerations.
\begin{Definition}\label{Def Final}
	We say that a linear manifold $\mL$ localizes the operator $A$ if the following holds:
	$$
	1)\,\,\,\mL\subset\Dom A,\quad 2)\,\,\,A\mL\subset\overline{\mL}.
	$$
	The operator $A_\mL:=A\harp\mL$ is called the part of $A$ on $\mL$.
\end{Definition}

Then a subspace $\mG$ is invariant for $A$ (according to Definition \ref{Def Our}) if the linear manifold $\mG\cap\Dom A$ is dense in $\mG$ and localizes $A$.

If $\mL\subset\mG$, $\overline{\mL}=\mG$, $\mL$ localizes $A$, and $\mG$ is invariant for $A$, then we obviously have $A_\mL\subset A_\mG$ (but not necessarily $A_\mL=A_\mG$). At the same time, a situation is possible where $\mL$ localizes $A$ but $A$ has no part $A_\mG$, i.e., there exist elements $x\in \mG\cap\Dom A$ such that $Ax\not\in\mG$. Let us give an example (A.\,F. Vakulenko).

Let $\mH:=L_2(\Bbb R_+)\oplus\Bbb C$, and define the operator $A_0$ as follows:
$$
A_0f:=\begin{pmatrix} -f''\\
	0
\end{pmatrix}; \qquad \Dom A_0:=\{f\in W^2_2(\Bbb R_+)\,|\,\,f(0)=0\}\oplus\Bbb C.
$$
Extend $A_0$ to $A$ on
$$
\Dom A:=\Dom A_0\dotplus \{\alpha\begin{pmatrix} e\\0
\end{pmatrix}\,|\,\,\alpha\in\Bbb C\},
$$
($e:=e^{-x}$, $x\in\Bbb R_+$), setting
$$
A\harp\Dom A_0:=A_0, \qquad A\begin{pmatrix} e\\0
\end{pmatrix}:=\begin{pmatrix} e\\1
\end{pmatrix}
$$
and extending $A$ to $\Dom A$ by linearity. The following facts are easily verified: 1) the operator $A$ has no part in the subspace $\mG:=L_2(\Bbb R_+)\oplus\Bbb C$; 2) the operator $A$ has a part on the linear manifold $\mL:= C^\infty_0(\Bbb R_+)\oplus\{0\}$, which is dense in $\mG$; 3) the linear manifold $\mL$ localizes $A$.
\medskip

\bul The issues discussed are quite relevant for applications. For example, an adequate understanding of the part of an operator is essential in controllability problems for dynamical systems with boundary control; see \cite{B_DSBC_3,BD_DSBC,BSim_ZNS_2019}.


\begin{thebibliography}{99}

\bibitem{AG_engl}
N.\,I. Akhiezer and I.\,M. Glazman.
\newblock {Theory of linear operators in Hilbert space. Volume I.}
\newblock {\it Dover Publishers, Inc. New-York}, 1992.

\bibitem{AG}
N.\,I. Akhiezer and I.\,M. Glazman.
\newblock {Theory of linear operators in Hilbert space.}
\newblock {\it Nauka, Moscow}, 1966. (in Russian)

\bibitem{Behrndt}
J. Behrndt, S.  Hassi, H. de Snoo.
\newblock{Boundary value problems, Weyl functions, and differential
	operators.}
\newblock{\it Birkh{\"a}user}, 2010; https://doi.org/10.1007/978-3-030-36714-5.

\bibitem{B_DSBC_3}
M.\,I. Belishev.
\newblock{Wave propagation in abstract dynamical
system with boundary control.}
\newblock{\it Zapiski Nauch. Semin. POMI}, 521 (2023), 8--32\,\,\,(in Russian).
{\it English translation}: arXiv:2307.00605v1 [math.DS] 2 Jul 2023.



\bibitem{BD_DSBC}
M.\,I. Belishev, M.\,N. Demchenko.
\newblock{Dynamical system with boundary control associated with
symmetric semibounded operator.}
\newblock{\it Journal of Mathematical Sciences},
October 2013, Volume 194 (2013), Issue 1, 8--20. DOI:
10.1007/s10958-013-1501-8.



\bibitem{BSim_ZNS_2019}
M.\,I. Belishev and S.\,A. Simonov. {\newblock On an evolutionary
dynamical system of the first soder with boundary control}.
\newblock{\it J. Math. Sci.}, (2021).
https://doi.org/10.1007/s10958-021-05183-y.

\bibitem{BSim_charact_SL}
M.\,I. Belishev, S.\,A. Simonov.
\newblock {A model and characterization of a class of symmetric semibounded
	operators.}
\newblock{\it arXiv:2504.01000, 1 Apr 2025}.



\bibitem{BirSol}
M.\,Sh. Birman, M.\,Z. Solomak.
\newblock{Spectral theory of self-adjoint operators in Hilbert space.}
\newblock{\it Reidel Publishing Comp.}, 1987.



\bibitem{DSch}
N. Dunford and J.\,T. Schwartz.
\newblock {Linear operators. Spectral theory. Self-adjoint operators in Hilbert space.}
\newblock {\it New-York, London}, 1963.



\bibitem{Ed}
R.\,E. Edwards.
\newblock{Functional analysis: theory and applications.}
\newblock{\it Dover Publications, Inc. Ney-York,} 1994.



\bibitem{KAk}
L.\,V. Kantorovich, G.\,P. Akilov.
\newblock{Functional analysis.}
\newblock{\it Nauka, Moscow}, 1977. (in Russian)


\bibitem{Kato}
T. Kato.
\newblock{Perturbation theory for linear operators.}
\newblock{\it Springer-Verlag, Berlin, Heidelberg, New-York}, 1966.



\bibitem{Naimark}
M.\,A. Naimark.
\newblock{Linear Differential Operators.}
\newblock{\it Harrap}, 1968.



\bibitem{Ples}
A.\,I. Plesner.
\newblock{Spectral Theory of Linear Operators.}
\newblock{\it Nauka, Moscow}, 1965. (in Russian)

\end{thebibliography}
\end{document}